\documentclass{amsart}
\usepackage[T1]{fontenc}
\usepackage[latin1]{inputenc}
\usepackage{hyperref}
\makeatletter

\theoremstyle{plain}    
\newtheorem{thm}{Theorem} 
\numberwithin{equation}{section} 
\numberwithin{figure}{section} 
\theoremstyle{plain}    
\newtheorem{Def}[thm]{Definition}
\newtheorem{lem}[thm]{Lemma} 
\theoremstyle{plain}    
\newtheorem{prop}[thm]{Proposition} 
\theoremstyle{remark}

\theoremstyle{remark}

\usepackage[all]{xy}
\xyoption{dvips}
\CompileMatrices

\newcommand{\cC}{\mathcal C}
\newcommand{\cT}{\mathcal T }

\newfont{\gothic}{eufm10 scaled 1100}
\newfont{\smgothic}{eufm10 scaled 900}

\font\tenCal=eusm10
\font\sevenCal=eusm7
\font\fiveCal=eusm5
\newfam\Calfam
  \textfont\Calfam=\tenCal
  \scriptfont\Calfam=\sevenCal
  \scriptscriptfont\Calfam=\fiveCal

\makeatother

\begin{document}

\title{Line Bundles on Complex Tori and a Conjecture of Kodaira}


\author{Jean-Pierre Demailly, Thomas Eckl, Thomas Peternell}



\address{Jean-Pierre Demailly, Institut Fourier, Universit\'e de Grenoble,
BP74 38402 St. Martin d'Heres, France}

\email{Jean-Pierre.Demailly@ujf-grenoble.fr}

\address{Thomas Eckl and Thomas Peternell, Institut für Mathematik, Universität Bayreuth, 
95440 Bayreuth, Germany}

\email{thomas.eckl@uni-bayreuth.de; thomas.peternell@uni-bayreuth.de}


\maketitle

\bibliographystyle{alpha}  

\section{A conjecture of Kodaira}
\label{Kod-section}

\noindent
A fundamental question in K\"ahler geometry asks whether any compact
K\"ahler manifold can be realised as a deformation of a projective
manifold. This is made more precise in the following

\begin{Def}
A compact K\"ahler manifold 
$X$ 
can be approximated algebraically or is almost algebraic if there exists a 
complex manifold 
$\mathcal{X}$ 
and a surjective holomorphic
submersion $\pi: \mathcal{X} \rightarrow \Delta$ 
to the unit disc 
$\Delta \subset \mathbb{C}$ 
such that 
$\mathcal{X}_0 \simeq X$ 
and there is a sequence 
$(t_k)$ 
converging to 
$0$ 
such that all 
$\mathcal{X}_{t_k}$ 
are projective. 
\end{Def}

\noindent
In \cite{Kod63} Kodaira proved that every K\"ahler surface is almost 
algebraic, and it is by now a standard conjecture, known as the Kodaira 
conjecture, that this should be true also in higher dimensions.

\noindent
One of the first things to try in higher dimensions is certainly to
look at projective bundles over K\"ahler non algebraic manifolds 
(possibly starting with surfaces for the sake of simplicity). 
Essentially the Kodaira conjecture implies that holomorphic vector 
bundles over certain K\"ahler manifolds should survive on sufficiently many 
algebraic approximations. The idea is based on the following easy
statement.

\begin{prop}
Let $X$ be a compact K\"ahler manifold which has a ${\mathbb{P}_r}$-bundle 
structure $X\to A$ over some complex torus $A$. Then for every 
deformation $\mathcal{X}\to S$ with $\mathcal{X}_0\simeq X$, the nearby 
fibers $\mathcal{X}_t$ have a ${\mathbb{P}_r}$-bundle structure
$\mathcal{X}_t\to \mathcal{A}_t$ where $\mathcal{A}$ is a deformation of 
$A$ in a neighborhood of $t=0$. Moreover, if $X=\mathbb{P}(V)$ for some 
vector bundle $V$ on $A$, then $\mathcal{X}_t=\mathbb{P}(V_t)$ for a 
suitable deformation $V_t\to \mathcal{A}_t$ of $V\to A$.
\end{prop}

\begin{proof} We look at the relative Albanese map $\alpha:\mathcal{X}\to
\mathcal{A}$. Then $\mathcal{A}\to S$ is a deformation of tori such that
$\alpha_t:\mathcal{X}_t\to\mathcal{A}_t$ is the Albanese map for each 
$t\in S$. Since $\alpha_0$ is a submersion, $\alpha_t$ should be also 
a submersion $t$ in a neighborhood $U\subset S$ of $0$, and the fibers 
of $\alpha_t$ are deformations of 
$\mathbb{P}_r$. Since $\mathbb{P}_r$ is undeformable, we conclude
that $\alpha_t:\mathcal{X}_t\to\mathcal{A}_t$  is also a $\mathbb{P}_r$-bundle
for small $t$. Now, the fact that $\mathcal{X}_t=\mathbb{P}(V_t)$ is equivalent
to the fact that the relative anticanonical bundle 
$K_{\mathcal{X}_t/\mathcal{A}_t}^{1}$ has an $(r+1)$-root $L_t$ on
$\mathcal{X}_t$, in which case $V_t=(\alpha_t)_*(L_t)$. However, the 
obstruction for a line bundle
to have an $(r+1)$-root lies in $H^2(\mathcal{X}_t, \mathbb{Z}/(r+1)
\mathbb{Z})$. This is a discrete locally constant coefficient system,
so if the obstruction vanishes for $t=0$, it must also vanish on the
connected component of $0$ in $U\subset S$.
\end{proof}

\noindent
Proposition 2 actually holds for arbitrary projective bundles over compact
K\"ahler manifolds; the proof is slightly more involved and is given in the
last section. Proposition 2 however is sufficient for our purposes.
\medskip 

\noindent
In view of this, it is
natural to look at the following potential candidate for a counter-example: 
Start with a $3$-dimensional complex torus with Picard number 
$\rho(A)  \geq 3$.
Let 
$L_i \in NS(A)$ 
be (numerical equivalence classes of) linearly independent holomorphic 
line bundles over 
$A$.
Let 
$U \subset \mathbb{C}^9$ 
be a neighborhood of 
$[A]$ 
in the universal deformation space of 
$A$.
As explained in the next section, every 
$L_i$ 
determines a 
$3$-codimensional subspace 
$V_i = V(L_i)$
in 
$U$ 
such that 
$c_1(L_i)$ 
is 
$(1,1)$, 
i.e. 
$L_i$ 
is a holomorphic line bundle on 
$A'$ 
if and only if 
$[A'] \in V_i.$ 
\medskip

\noindent
{\it Now we make the following Assumption:\\
The intersection of the
$V_i$'s
has the expected dimension 0, i.e.
\[ (*)\kern30pt V_1 \cap V_2 \cap V_3~~
\mbox{contains $\{A\}$ as an isolated point}.\] }
Then consider the 
$6$-dimensional manifold
\[ Y = \mathbb{P}(\mathcal{O}_A \oplus L_1) \times_A 
       \mathbb{P}(\mathcal{O}_A \oplus L_2) \times_A 
       \mathbb{P}(\mathcal{O}_A \oplus L_3). \]
This is a 
$\mathbb{P}_1^3$-bundle over 
$A$ 
with projection 
$\pi: Y \rightarrow A.$
In each subspace 
$\mathbb{P}(\mathcal{O}_A \oplus L_i)$ 
there is a section 
$Z_i$ 
at infinity given by the direct summand 
$\mathcal{O}_A$.
This gives a section 
$Z$ 
of 
$\pi$ 
by selecting over every 
$a \in A$
the point 
$(x_1,x_2,x_3)$, 
where 
$\{x_i\} = Z_i \cap \pi^{-1}(a)$.

\begin{prop}
The blow up
$\sigma: X \rightarrow Y$
of
$Z \subset Y$
is rigid in the sense that there is no positive-dimensional family of 
deformations of 
$X$.
\end{prop}
\begin{proof}
Notice that, denoting by
$\mathbb{P}_1^3(x)$ 
the blow up of
$\mathbb{P}_1^3$ 
in one point, 
$X$ 
is a 
$\mathbb{P}_1^3(x)$-bundle over 
$A$. 
So let 
$(X_t)$ 
be a deformation of 
$X = X_0$ 
over the 
$1$-dimensional unit disc 
$\Delta$.
The first step is to proof that, possibly after shrinking 
$\Delta$, 
every 
$X_t$ 
is a 
$\mathbb{P}_1^3(x)$-bundle over its 
($3$-dimensional) Albanese torus 
$A_t$. 
In fact, 
$q(X_t) = 3$ 
for all 
$t$ 
and the Albanese map 
$\alpha_t$ 
is smooth for small 
$t$. 
In order to prove that 
$\alpha_t$ 
is a 
$\mathbb{P}_1^3(x)-$
bundle, it suffices to show that 
$\mathbb{P}_1^ 3(x)$ 
is rigid, i.e. every small deformation of
$\mathbb{P}_1^3(x)$ 
is again 
$\mathbb{P}_1^3(x)$. 

\noindent
In fact, let 
$Z = \mathbb{P}_1^3(x)$ 
for simplicity of notations. Let 
$\tau: Z \rightarrow \mathbb{P}_1^3$ 
be the blow-up map with exceptional divisor 
$E \simeq \mathbb{P}_2$. 
Then there is an exact sequence
\[ 0 \rightarrow T_Z \rightarrow \tau^*T_{\mathbb{P}_1^3} 
     \rightarrow T_E(-1) \rightarrow 0. \]
Since 
$\dim H^0(T_{\mathbb{P}_1^3}) = 9$, 
$\dim H^0(T_Z) = 6$,
$H^0(T_E(-1)) = 3$ 
and 
$H^1(\tau^*T_{\mathbb{P}_1^3}) = 0$,
by taking cohomology of the above exact sequence it follows
\[ H^1(T_Z) = 0, \]
in particular 
$Z$ 
is rigid.

\noindent
Let 
$\mathcal{X}$ 
be the total space of 
$(X_t)$ 
and let 
$\pi: \mathcal{X} \rightarrow \mathcal{A}$ 
be the relative Albanese
map for 
$\mathcal{X} \rightarrow \Delta$. 
Then 
$\mathcal{A} \rightarrow \Delta$ 
is a torus bundle; let 
$A_t$
be the fiber over
$t$, 
so that 
$A = A_0.$ 
Now the exceptional divisor 
$D$ 
of 
$\sigma$ 
moves in 
$\mathcal{X}$. 
This is easy to see by considering 
$D \cap \pi^{-1}(a) = \mathbb{P}_2$ 
for 
$a \in A$. 
In fact, the normal bundle of this 
$\mathbb{P}_2 $ 
is 
$\mathcal{O}(-1) \oplus \mathcal{O}$, 
so that the 
$\mathbb{P}_2$ 
moves and forces
$D$ 
to move. Therefore one obtains a fiberwise blow-down 
$\mathcal{X} \rightarrow \mathcal{Y}$ 
inducing the birational map
$\sigma: X \rightarrow Y$.
Of course there is a factorisation 
$\mathcal{X} \rightarrow \mathcal{Y} \rightarrow \mathcal{A}$ 
and 
$\mathcal{Y} \rightarrow \mathcal{A}$ 
is a 
$\mathbb{P}_1^3$-bundle. Again let 
$Y_t$ 
be the fiber over 
$t$.
Next it is shown that it is possible to write
\[  Y_t = Y_{1,t} \times_{A_t} Y_{2,t} \times_{A_t} Y_{3,t} \]
with 
$\mathbb{P}_1$-bundles 
$Y_{i,t} / A_t$, 
and this can be done simultaneously, i.e. the
$(Y_{i,t})$ 
form a family 
$\mathcal{Y}_i$. 
The most economic way to do that is to note that the relative 
Picard number 
$\rho(\mathcal{Y}/\mathcal{A})$ 
equals
$3$
since 
$\rho(Y_0/A_0) = 3$ 
(this is a product situation).
By taking relative extremal contractions in the sense of Mori theory one gets 
a tower of three 
$\mathbb{P}_1$-bundles. Of course there are three choices of the first one and then two 
choices for the second since the situation is completely symmetric in 
$i$. (This situation could possibly lead to some monodromy action
$\pi_1(A_t)\to \frak{S}_3$, but since such actions are discrete and
depend continuoulsly on $t$, the fact that we have a non twisted
product for $t=0$ implies that we have no such twist for $t$ arbitrary). 
The last 
contraction will provide the space 
$\mathcal{Y}_i$ 
for the appropriate 
$i$. 
Now consider the canonical map 
\[ Y_t \rightarrow Y_{1,t} \times_{A_t} Y_{2,t} \times_{A_t} Y_{3,t}. \]
Then this map is immediately seen to be an isomorphism. 

\noindent
Since 
$Y_{i,t}$ 
is a 
$\mathbb{P}_1$-bundle over 
$A_t$ 
and since it is has a section by construction, it follows 
\[ Y_{i,t} = \mathbb{P}(E_{i,t}) \]
with a rank 2-bundle 
$E_{i,t}$
(normalized such that 
$E_{0,t} = \mathcal{O}_{A_0} \oplus L_i$), 
and the 
$E_{i,t}$ 
form a holomorphic rank 2-bundle 
$\mathcal{E}_i$
over 
$\mathcal{A}$. 
Since the section at infinity in 
$Y_0$ 
deforms by construction to sections in 
$Y_t$, 
one obtains a global quotient 
$\mathcal{E}_i \rightarrow \mathcal{G}_i \rightarrow 0$ 
such that 
$\mathcal{G}_i | A_0 = \mathcal{O}_{A_0}$. 
By changing
$\mathcal{E}_i$ 
appropriately, one may assume that 
$\mathcal{G}_i = \mathcal{O}_{\mathcal{A}}$. 
Let 
$\mathcal{L}_i$ 
be the kernel of 
$\mathcal{E}_i \rightarrow \mathcal{O}_{\mathcal{A}}$.
Then 
$\mathcal{L}_i | A_0 = L_i$. 
But this implies that there is a deformation of 
$A = A_0$
such that all three line bundles 
$L_i$ 
remain holomorphic. But the assumption 
\[ V_1 \cap V_2 \cap V_3 = \{A\} \]
implies that there is no such (nontrivial) deformation of
$A$. 
\end{proof}

\noindent
It is therefore a very natural question to ask whether these rigid 
$6$-dimensional K\"ahler manifolds are projective or not. If they were
not projective, we would get counter-examples to the Kodaira conjecture.
Unfortunately (in view of getting easy counter-examples!), Theorem~3 
of the next section tells us that a complex torus $A$ verifying Assumption 
(*) for some triple of holomorphic line bundles $L_i$ is always
an abelian variety.

\section{Holomorphic line bundles on complex tori}

\noindent
Let
$X$
be a complex torus of dimension 
$g$. 
As explained in
\cite{CT}, \cite{CAV}
$X$
admits a period matrix of the form 
$(\tau, \mathbf{1}_g)$
with
$\tau \in M_g(\mathbb{C})$,
the 
$g \times g$-matrices with entries in
$\mathbb{C}$
such that
$\mathrm{det}(\mathrm{Im}\ \tau) \neq 0$.
Conversely every such matrix is the period matrix of a complex torus.

\noindent
If 
$\Lambda \in \mathbb{C}^g := V$
denotes the lattice generated by the columns of
$(\tau, \mathbf{1}_g)$
the N\'eron-Severi group of
$X$
may be described as
\[ NS(X) = \left\{ E = \left( \begin{array}{cc}
                          A & B \\
                          - {^t}B & C               
                          \end{array} \right)
                   \in M_{2g}(\mathbb{Z})
                   \left| 
                   \begin{array}{l}
                   A\ \mathrm{and}\ C\ \mathrm{alternating,\ and\ } \\ 
                   A - B\tau + {^t}\tau {^t}B + {^t}\tau C \tau = 0 
                   \end{array} \right. \right\}. \]
The equality ensures that the alternating form 
$E$
is a 
$(1,1)$-form, cf.\ \cite[p.$\,$10]{CT}.
\\
\noindent
\begin{thm} \label{alg-thm}
Let
$X$
be a 3-dimensional complex torus with period matrix
$(\tau,\mathbf{1}_3)$
and let
$E_1 \cdot \mathbb{Z} \oplus E_2 \cdot \mathbb{Z} \oplus E_3 \cdot 
                                          \mathbb{Z}\ \subset NS(X)$
be a rank 3 subgroup of the N\'eron-Severi group 
$NS(X)$
of
$X$
generated by
$E_1, E_2, E_3 \in NS(X)$.
Then there is a sequence
$(X_n)$
of 3-dimensional complex tori with period matrices
$(\tau_n, \mathbf{1}_3)$
such that
\begin{itemize}
\item[(i)]
the 
$\tau_n$
converge to
$\tau$
for 
$n \rightarrow \infty$,
\item[(ii)]
$ E_1 \cdot \mathbb{Z} \oplus E_2 \cdot \mathbb{Z} \oplus E_3 \cdot 
                                          \mathbb{Z} \subset NS(X_n)$
and
\item[(iii)]
$X_n$
is a complex abelian variety.
\end{itemize}
\end{thm}

\noindent
As a first step towards a proof,
$E = \left( \begin{array}{cc}
            A & B \\
      - {^t}B & C               
            \end{array} \right)$
may be considered as an element of the free abelian group 
$\mathbb{Z}^{15}$: 
the matrices
$A = \left( \begin{array}{ccc}
            0 & a_1 & a_2 \\
         -a_1 & 0   & a_3 \\
         -a_2 & -a_3 & 0         
            \end{array} \right)$
and
$C = \left( \begin{array}{ccc}
            0 & c_1 & c_2 \\
         -c_1 & 0   & c_3 \\
         -c_2 & -c_3 & 0         
            \end{array} \right)$
are alternating, and
$B = \left( \begin{array}{ccc}
            b_1 & b_2 & b_3 \\
            b_4 & b_5 & b_6 \\
            b_7 & b_8 & b_9         
            \end{array} \right)$
is arbitrary. Since
$k \cdot E \in NS(X)$
implies
$E \in NS(X)$,
condition (ii)
is equivalent to
\[ E_1 \cdot \mathbb{Q} \oplus E_2 \cdot \mathbb{Q} \oplus E_3 \cdot 
                                          \mathbb{Q} \subset 
   NS(X_n) \otimes_\mathbb{Z} \mathbb{Q}, \]
and
$E_1 \cdot \mathbb{Q} \oplus E_2 \cdot \mathbb{Q} \oplus E_3 \cdot \mathbb{Q}$
may be interpreted as a 
$\mathbb{Q}$-rational point in the Grassmannian
$G(3,15)$.

\noindent
For a given 3-dimensional subspace
$E_1 \cdot \mathbb{Q} \oplus E_2 \cdot \mathbb{Q} \oplus E_3 \cdot 
                                          \mathbb{Q} \subset \mathbb{Q}^{15}$
the equations 
$A_i - B_i\tau + {^t}\tau {^t}B_i + {^t}\tau C_i \tau = 0$,
$i = 1,2,3$
imply algebraic relations between the entries of 
\[ \tau = \left( \begin{array}{ccc}
            \tau_1 & \tau_2 & \tau_3 \\ 
            \tau_4 & \tau_5 & \tau_6 \\ 
            \tau_7 & \tau_8 & \tau_9          
               \end{array} \right): \]
Since the
$A_i - B_i\tau + {^t}\tau {^t}B_i + {^t}\tau C_i \tau$
are alternating matrices, the number of these  relations can be reduced to 9
($i = 1,2,3$):
\[ (\ast)\ \ 
   \begin{array}{rcl} 
   0& = & a_{i1} - b_{i1}\tau_2 - b_{i2}\tau_5 - b_{i3}\tau_8 + b_{i4}\tau_1 + 
            b_{i5}\tau_4 + b_{i6}\tau_7 \\ 
    & & \mbox{} + c_{i1}(\tau_1\tau_5 - \tau_2\tau_4) +  
                  c_{i2}(\tau_1\tau_8-\tau_2\tau_7) +
                  c_{i3}(\tau_4\tau_8-\tau_5\tau_7)  \\
   0& = & a_{i2} - b_{i1}\tau_3 - b_{i2}\tau_6 - b_{i3}\tau_9 + b_{i7}\tau_1 + 
            b_{i8}\tau_4 + b_{i9}\tau_7 \\ 
    & & \mbox{} + c_{i1}(\tau_1\tau_6 - \tau_3\tau_4) +  
                  c_{i2}(\tau_1\tau_9-\tau_3\tau_7) +
                  c_{i3}(\tau_4\tau_9-\tau_6\tau_7)  \\
   0& = & a_{i3} - b_{i4}\tau_3 - b_{i5}\tau_6 - b_{i6}\tau_9 + b_{i7}\tau_2 + 
            b_{i8}\tau_5 + b_{i9}\tau_8 \\ 
    & & \mbox{} + c_{i1}(\tau_2\tau_6 - \tau_3\tau_5) +  
                  c_{i2}(\tau_2\tau_9-\tau_3\tau_8) +
                  c_{i3}(\tau_5\tau_9-\tau_6\tau_8).  \\
   \end{array} \] 
So there is an algebraic subset 
$U_{E_1,E_2,E_3}$
of
$\mathbb{C}^9 = \mathbb{C}^3 \times \mathbb{C}^3$
such that
$U_{E_1,E_2,E_3} \cap \left\{ \tau \in \mathbb{C}^9 : 
                              \mathrm{det}(\mathrm{Im\ }\tau) \neq 0 \right\}$
describes all
$\tau$'s
with
$E_1 \cdot \mathbb{Q} \oplus E_2 \cdot \mathbb{Q} \oplus E_3 \cdot 
                                          \mathbb{Q} \subset 
        NS(X_\tau) \otimes_\mathbb{Z} \mathbb{Q}$
where
$X_\tau$
is the complex torus corresponding to the period matrix
$(\tau, \mathbf{1}_3)$.
In particular, the union of all these
$U_{E_1,E_2,E_3}$
is an algebraic family
$U \subset G(3,15) \times \mathbb{C}^9$.
Let
$\bar{U} \subset G(3,15) \times \mathbb{P}^9$
denote the projective closure of
$U$.

\noindent
The heart of the proof is now a careful analysis of this family
$\bar{U}$,
especially of the fibers over
$\mathbb{Q}$-rational points of
$G(3,15)$. 
If they always contain an (analytically) dense subset of
$\tau$'s 
such that
$X_\tau$
is a complex abelian variety, the theorem will follow.

\noindent
The first observation is that all coefficients in the equations of
$(\ast)$
are rational. Hence,
$\mathbb{Q}$
is the field of definition of
$\bar{U}$, 
i.e. there exists a
$\mathbb{Q}$-scheme
$\bar{U}_\mathbb{Q}$
such that
$\bar{U} = \bar{U}_\mathbb{Q} \times_\mathbb{Q} \mathrm{Spec\ }\mathbb{C}$.
In particular, every fiber of
$\bar{U}$
over a
$\mathbb{Q}$-rational point of 
$G(3,15)$
has 
$\mathbb{Q}$
as field of definition, too.

\noindent
Next, one computes a fiber
$\bar{U}_{E_1,E_2,E_3}$
of
$\bar{U}$
with sufficiently general entries in the matrices
$E_1,E_2,E_3$.
This can be done with the computer algebra program Macaulay2 (\cite{M2},
\cite{AGM2}). Setting
\[ A_1 = \left( \begin{array}{ccc}
         0 & 0 & 0 \\
         0 & 0 & 2 \\
         0 & -2 & 0
         \end{array} \right) ,\  
   B_1 = \left( \begin{array}{ccc}
         1 & 1 & 0 \\
         1 & 1 & 2 \\
         1 & 1 & 2
         \end{array} \right) ,\  
   C_1 = \left( \begin{array}{ccc}
         0 & 1 & 0 \\
         -1 & 0 & 0 \\
         0 & 0 & 0
         \end{array} \right) , \]

\[ A_2 = \left( \begin{array}{ccc}
         0 & 1 & 2 \\
         -1 & 0 & 1 \\
         -2 & -1 & 0
         \end{array} \right) ,\  
   B_2 = \left( \begin{array}{ccc}
         0 & 0 & 0 \\
         1 & 1 & 1 \\
         0 & 1 & 0
         \end{array} \right) ,\  
   C_2 = \left( \begin{array}{ccc}
         0 & 0 & 0 \\
         0 & 0 & 1 \\
         0 & -1 & 0
         \end{array} \right) , \]

\[ A_3 = \left( \begin{array}{ccc}
         0 & 1 & 2 \\
         -1 & 0 & 1 \\
         -2 & -1 & 0
         \end{array} \right) ,\  
   B_3 = \left( \begin{array}{ccc}
         1 & 1 & 1 \\
         1 & 2 & 1 \\
         1 & 2 & 1
         \end{array} \right) ,\  
   C_3 = \left( \begin{array}{ccc}
         0 & 0 & 1 \\
         0 & 0 & 0 \\
         -1 & 0 & 0
         \end{array} \right) , \]

\noindent
(the matrix entries were chosen by a random number generator) and using the 
following Macaulay2 script
\begin{verbatim}
   k = QQ;  
   PT = k[t_0..t_9];

   A1 = matrix(PT,{{0, 0, 0}, {0, 0, 2}, {0, -2, 0}});
   B1 = matrix(PT,{{1,1,0},{1,1,2},{1,1,2}});
   C1 = matrix(PT,{{0,1,0},{-1,0,0},{0,0,0}});

   A2 = matrix(PT,{{0,1,2},{-1,0,1},{-2,-1,0}});
   B2 = matrix(PT,{{0,0,0},{1,1,1},{0,1,0}}); 
   C2 = matrix(PT,{{0,0,0},{0,0,1},{0,-1,0}});

   A3 = matrix(PT,{{0,1,2},{-1,0,1},{-2,-1,0}});
   B3 = matrix(PT,{{1,1,1},{1,2,1},{1,2,1}});
   C3 = matrix(PT,{{0,0,1},{0,0,0},{-1,0,0}});

   gent = genericMatrix(PT,t_1,3,3);

   s1 = matrix(PT,{{t_0,0,0},{0,t_0,0},{0,0,t_0}});
   s2 = s1*s1;

   Q1 = A1*s2 - B1*gent*s1 + transpose(gent)*transpose(B1)*s1 +\
   transpose(gent)*C1*gent;

   Q2 = A2*s2 - B2*gent*s1 + transpose(gent)*transpose(B2)*s1 +\
   transpose(gent)*C2*gent;

   Q3 = A3*s2 - B3*gent*s1 + transpose(gent)*transpose(B3)*s1 +\
   transpose(gent)*C3*gent;
   Q = Q1|Q2|Q3;
   --- Q contains the 9 relations between the t_i's homogenized with \
   respect to t_0

   q = saturate(ideal(flatten Q), ideal(t_0))
   --- saturation with t_0 removes all components on the hyperplane\
   t_0 = 0
   betti q
\end{verbatim}
one gets 8 linear and 1 quadratic equation describing the projective closure of
$U_{E_1,E_2,E_3}$: 
\begin{verbatim}
   t_7+3/5t_8+8/5t_9        
   t_6-3/20t_8+1/10t_9       
   t_5-3/5t_8+2/5t_9         
   t_4+1/2t_8+t_9            
   t_3-1/20t_8-3/10t_9       
   t_2+3/10t_8+9/5t_9        
   t_1-3/10t_8+1/5t_9        
   t_0-1/4t_8-3/2t_9         
   t_8^2-48t_8t_9-172/3t_9^2 
\end{verbatim}
Since the quadratic generator has discriminant 
$24^2+4\frac{172}{3}>0$
which is not the square of a rational number, this is a 
$\mathbb{Q}$-irreducible 0-dimensional scheme of degree 2; over
$\mathbb{C}$
it consists of two points.

\noindent
Unfortunately, these equations may cut out too much, since the projective 
closure of a fiber may be less than the fiber of the projective closure of a 
family. To deal with this problem one has to do a little detour: First one 
looks at the (inhomogeneous) ideal of the whole family 
$U$:
\begin{verbatim}
   k = QQ;
   P = k[t_0..t_9];
   PE = k[e_0..e_11,f_0..f_11,g_0..g_11];
   PT = P ** PE;

   A1 = matrix(PT,{{0, e_0, e_1}, {-e_0, 0, e_2}, {-e_1, -e_2, 0}});
   B1 = matrix(PT,{{e_3,e_4,e_5},{e_6,e_7,e_8},{e_9,e_10,e_11}});
   C1 = matrix(PT,{{0,1,0},{-1,0,0},{0,0,0}});

   A2 = matrix(PT,{{0, f_0, f_1}, {-f_0, 0, f_2}, {-f_1, -f_2, 0}});
   B2 = matrix(PT,{{f_3,f_4,f_5},{f_6,f_7,f_8},{f_9,f_10,f_11}});
   C2 = matrix(PT,{{0,0,0},{0,0,1},{0,-1,0}});

   A3 = matrix(PT,{{0, g_0, g_1}, {-g_0, 0, g_2}, {-g_1, -g_2, 0}});
   B3 = matrix(PT,{{g_3,g_4,g_5},{g_6,g_7,g_8},{g_9,g_10,g_11}});
   C3 = matrix(PT,{{0,0,1},{0,0,0},{-1,0,0}});

   gent = genericMatrix(PT,t_1,3,3);

   Q1 = A1 - B1*gent + transpose(gent)*transpose(B1) +\
   transpose(gent)*C1*gent;
   Q2 = A2 - B2*gent + transpose(gent)*transpose(B2) +\   
   transpose(gent)*C2*gent;
   Q3 = A3 - B3*gent + transpose(gent)*transpose(B3) +\
   transpose(gent)*C3*gent;

   Q = Q1|Q2|Q3;
   q = ideal flatten Q;
\end{verbatim}

\noindent
The projective closure of 
$U$
may be determined by computing a Groebner basis of this ideal with respect to
a monomial order refining the order by degree in the 
$t_i$'s 
and then homogenizing the generators with respect to
$t_0$ 
(\cite[15.31]{ECA}). This computation is too complicated for the whole 
Groebner basis, but it is already enough to look at the first few elements 
which are added to the original generators:
\begin{verbatim}
   gbasis = gb(q,PairLimit=>31);

   hgbasis = homogenize(gens gbasis,t_0,{1,1,1,1,1,1,1,1,1,1, 0,0,0,\
   0,0,0,0,0,0,0,0,0, 0,0,0,0,0,0,0,0,0,0,0,0, 0,0,0,0,0,0,0,0,0,0,\ 
   0,0});
\end{verbatim}  
Evaluation at 
$(E_1,E_2,E_3)$
\begin{verbatim}
   f = map(PT,PT,matrix(PT, {{t_0,t_1,t_2,t_3,t_4,t_5,t_6,t_7,t_8,\
   t_9, 0,0,2,1,1,0,1,1,2,1,1,2, 1,2,1,0,0,0,1,1,1,0,1,0, 1,2,1,1,\
   1,1,1,2,1,1,2,1}}));

   genfibre = ideal f(hgbasis);

   betti gb genfibre
\end{verbatim}
shows that the fibre
$(\overline{U})_{E_1,E_2,E_3}$
is contained in a scheme cut out by 8 linear and 1 quadratic equation, so 
\[ (\overline{U})_{E_1,E_2,E_3} = U_{E_1,E_2,E_3}. \]

\noindent
One can get further informations about
$\overline{U}$
from the homogenized equations collected in 
\texttt{hgbasis}. 
Since the projective closure of a fiber is equal to the fiber of the 
projective closure on an open subset they contain 9 equations describing the 
fibers of
$\overline{U}$
over an open subset around
$(E_1,E_2,E_3)$. 
Furthermore the command
\begin{verbatim}
   transpose leadTerm hgbasis
\end{verbatim}
shows that all of these equations contain 
$t$-variables. Hence each of these fibers is cut out by 9 non-constant equations, 
so it is not empty. Consequently, 
\textbf{no}
fiber is empty.

\noindent
Turning to the fibers of
$\overline{U}$
over
$\mathbb{C}^9 = \mathbb{C}^3 \times \mathbb{C}^3$
(resp. 
$\mathbb{P}^9$) 
one sees immediately that these are non-empty linear subspaces. Hence
$\overline{U}$
is connected. Finally,
$U$
is regular, as may be shown by deriving the equations in 
$(\ast)$
with respect to the
$a_{ij}$'s.
Taken all these facts together it follows that
$U$
and hence
$\overline{U}$
is irreducible. So by Stein factorization every 
$0$-dimensional fiber is of degree 2.

\noindent
Now it is easy to prove for these 0-dimensional fibers over
$\mathbb{Q}$-rational points that they describe period matrices
$\tau$
belonging to complex abelian varieties: Since the fibers are
$\mathbb{Q}$-rational, too, the entries of
$\tau$
are elements of a field extension of
$\mathbb{Q}$
of degree 2. The defining equations of the N\'eron-Severi group show that then
$NS(X_\tau)$
is a
$15 - 2 \times 3 = 9$-dimensional
$\mathbb{Q}$-vector space. But a 3-dimensional complex torus with maximal Picard number 9 is
algebraic (cf. \cite{CT}).

\noindent
What about the higher dimensional fibers ? Consider the 
$\mathbb{Q}$-rational map
$\phi: G(3,15) \dasharrow \mathrm{Hilb}^2(\mathbb{P}^9_\mathbb{Q})$
whose existence is the essence of the arguments used above. Let
\[ \xymatrix{
     G \ar[d]_{\pi} \ar[dr]^{\overline{\phi}} &  \\
     G(3,15) \ar@{->}[r]^{\phi} & {\mathrm{Hilb}^2(\mathbb{P}^9_\mathbb{Q})}  
            } \]
be the resolution of the singularities of
$\phi$
by blowing up
$\mathbb{Q}$-regular centers. This is possible by the Hironaka package, see 
\cite{Hir64}, \cite{BM97} or \cite{RS97}. Now the theorem is a consequence 
of the following
\\
\begin{lem}
Let 
$Z \subset Y$
be an embedding of
$\mathbb{Q}$-regular 
$\mathbb{Q}$-schemes, let 
$z \in Z$
be a 
$\mathbb{Q}$-rational point and let
$\phi: \tilde{Y} \rightarrow Y$
be the blow up of 
$Y$
with center
$Z$. Then the 
$\mathbb{Q}$-rational points are dense on the fiber
$\phi^{-1}(z)$.
\end{lem}
\begin{proof}
This is almost trivial: Choose a
$\mathbb{Q}$-regular sequence
$(f_1, \ldots, f_s,f_{s+1}, \ldots, f_t)$
in the local ring
$\mathcal{O}_{Y,z}$
such that
$\hbox{\gothic m}_{Z,z} = (\overline{f}_{s+1}, \ldots, \overline{f}_t )
 \subset \mathcal{O}_{Z,z}$
and
$\hbox{\gothic m}_{Y,z} = (f_1, \ldots, f_t )$.
The blowing up of 
$\mathrm{Spec} \mathcal{O}_{Y,z}$
with center
$\mathrm{Spec} \mathcal{O}_{Z,z}$
is given by
\[ \mathrm{Proj}\ \mathcal{O}_{Y,z}[f_{s+1}, \ldots, f_t] =
   (\mathrm{Spec}\ \mathcal{O}_{Y,z} \times \mathbb{P}_\mathbb{Q}^{t-s-1})/
   V(T_if_j - T_jf_i), \]
and the fiber over z is
$\cong \mathbb{P}_\mathbb{Q}^{t-s-1}$.
\end{proof}

\noindent
Apply the lemma on
$\phi$:
If
$[E_1 \cdot \mathbb{Q} \oplus E_2 \cdot \mathbb{Q} \oplus 
  E_3 \cdot \mathbb{Q}] = [W] \in G(3,15)$
is a 
$\mathbb{Q}$-rational point then
$\pi^{-1}([W]) \subset G$
will contain an analytically dense subset of
$\mathbb{Q}$-rational points, and the same will be true of the image
$\overline{\phi}(\pi^{-1}([W])) \subset 
 \mathrm{Hilb}^2(\mathbb{P}^9_\mathbb{Q})$.
But 
$\mathbb{Q}$-rational points in
$\mathrm{Hilb}^2(\mathbb{P}^9_\mathbb{Q})$
describe pairs of points corresponding to abelian varieties, and all pairs in
$\overline{\phi}(\pi^{-1}([W]))$
map surjectively on the fiber over 
$[W]$
in
$\overline{U}$.
Hence this fiber contains a dense open subset of period matrices 
$\tau$
such that 
$X_\tau$
is an abelian variety.
\\
\noindent
\textbf{Remark.}
Some words about the Macaulay2 computations: Since all the relevant equations 
and varieties are defined over
$\mathbb{Q}$
and also the operations applied to them like taking the projective closure 
work over
$\mathbb{Q}$,
these calculations give exact results.

\section{Modifications and a general setting for counter-examples}

\noindent
Of course the construction in section~\ref{Kod-section} possibly could be modified in 
several ways and then might lead to a counter-example to Kodaira's conjecture.

\noindent
First we show that in the setup before Proposition 2 - of course now without Assumption (*) -
the variety $X$ from Proposition 2 can be approximated algebraically. In deed in that situation (using the
old terminology),
$V_1 \cap V_2 \cap V_3$
contains other complex tori than
$A$.
Then theorem~\ref{alg-thm} assures the existence of a sequence 
$\left\{ A_n \right\}_{n \in \mathbb{N}} \subset V_1 \cap V_2 \cap V_3$
of 
abelian varieties converging to 
$A$.
The following lemma shows that this implies
$X$
almost algebraic:
\begin{lem}
Let
$E = \left( \begin{array}{cc}
                          A & B \\
                          - {^t}B & C               
                          \end{array} \right) \in M_{2g}(\mathbb{Z})$
be a skew symmetric matrix with integral entries and let
\[ V = \left\{ \tau \in M_3(\mathbb{C}) |  
                     A - B\tau + {^t}\tau {^t}B + {^t}\tau C \tau = 0;\
               \det \mathrm{Im}\tau \neq 0 \right\}
             \subset \mathbb{C}^9 \]
be the set of period matrices 
$\tau$
such that
$X_\tau$
is a complex torus with
$E \in NS(X_\tau)$.

\noindent
Let 
$\mathcal{X} = V \times \mathbb{C}^g/\Lambda_\tau$
be the family of these tori
$X_\tau$ 
where
$\Lambda_\tau = (\tau, \mathbf{1}_g)$
is the lattice belonging to
$X_\tau = \mathbb{C}^g/\Lambda_\tau$.
Then every
$\tau_0 \in V$
has an open neighborhood
$U \subset V$
such that there exists a holomorphic line bundle
$\mathcal{L}_U$
on
$\mathcal{X}_U$
such that 
$c_1(\mathcal{L}_\tau) = E$
for all
$\tau \in U$.
\end{lem}
\begin{proof}
Let
$\pi: \mathcal{X} \rightarrow V$
be the projection of
$\mathcal{X}$
onto
$V$.
By taking direct images with respect to
$\pi$
and deriving the long exact sequence from
$0 \rightarrow \mathbb{Z} \rightarrow \mathcal{O}_{\mathcal{X}} \rightarrow   
    \mathcal{O}_{\mathcal{X}}^\ast \rightarrow 1$
one obtains the sequence
\[ R^1\pi_\ast \mathcal{O}_{\mathcal{X}}^\ast \rightarrow 
   R^2\pi_\ast \mathbb{Z} \rightarrow R^2\pi_\ast \mathcal{O}_{\mathcal{X}}. \]
The skew symmetric matrix
$E$
gives a section of
$R^2\pi_\ast \mathbb{Z}$
which is mapped to 0 in
$R^2\pi_\ast \mathcal{O}_{\mathcal{X}}$
since
$E \in NS(X_\tau)$
for all
$\tau \in V$.
Hence 
$E$
is the image of a section in
$R^1\pi_\ast \mathcal{O}_{\mathcal{X}}^\ast$. 
Take an open neighborhood 
$U$
of
$\tau_0$
such that the section restricted to 
$U$
is a cohomology class in 
$H^1(\mathcal{X}_{|\pi^{-1}(U)}, \mathcal{O}_{\mathcal{X}}^\ast)$.
This class gives the line bundle 
$\mathcal{L}_U$.
\end{proof}

\noindent
Next, consider the following more 
general setting: Take an 
$n$-dimensonal complex torus 
$A$
and
$k$
vector bundles
$E_1, \ldots , E_k$
over
$A$
of rank 
$r_1, \ldots , r_k \leq n$.
Let 
$Y$
be the
$(n + r_1 + \ldots + r_k)$-dimensional manifold
\[ \mathbb{P}(\mathcal{O}_A \oplus E_1) \times_A \cdots \times_A 
   \mathbb{P}(\mathcal{O}_A \oplus E_k). \]
This a 
$(\mathbb{P}^{r_1} \times \ldots \times \mathbb{P}^{r_k})$-bundle over
$A$
with projection
$\pi: Y \rightarrow A$.
In each subspace 
$\mathbb{P}(\mathcal{O}_A \oplus E_i)$ 
there is a section 
$Z_i$ 
at infinity given by the direct summand 
$\mathcal{O}_A$.
This gives a section 
$Z$ 
of 
$\pi$ 
by selecting over every 
$a \in A$
the point 
$(x_1,\ldots,x_k)$, 
where 
$\{x_i\} = Z_i \cap \pi^{-1}(a)$.
Let
$\sigma: X \rightarrow Y$
be the blow up of
$Z \subset Y$.
\\
\noindent
Similar arguments as in section~\ref{Kod-section} show
\begin{prop}
If there is a positive-dimensional family of deformations of 
$X$
then there will also exist a deformation family of complex tori
$\left\{ A_t \right\}_{t \in \Delta}$
such that 
$A = A_0$
and all vector bundles 
$E_1, \ldots , E_k$
remain holomorphic on 
$A_t$. 
\hfill $\Box$
\end{prop}

\noindent
The condition on the vector bundles to remain holomorphic requires some further explanations: Let
$E$
be a vector bundle of rank
$r$
over an 
$g$-dimensional torus
$A$.
Then the Chern classes
$c_i(E)$
are
$(i,i)$-classes in
$H^{2i}(A,\mathbb{Z}) = \bigwedge^{2i}\mathrm{Hom}(\Lambda, \mathbb{Z})$,
where 
$\Lambda \subset \mathbb{C}^g =: V$
is a (non-degenerate) lattice such that 
$A = V/\Lambda$. 
Since
$H^{i,i}(A,\mathbb{C}) = 
 \bigwedge^i\mathrm{Hom}_\mathbb{C}(V,\mathbb{C}) \times
 \bigwedge^i\mathrm{Hom}_{\overline{\mathbb{C}}}(V,\mathbb{C}),$
the
$(i,i)$-classes in
$H^{2i}(A,\mathbb{Z})$
may be interpreted as a real valued alternating form
$F$
on 
$\bigwedge^i V$
such that
\[ F(i\Phi, i\Psi) = F(\Phi, \Psi)\ \mathrm{and}\ 
   F(\bigwedge^i \Lambda, \bigwedge^i \Lambda) \subset \mathbb{Z}. \]

\noindent
As in the case of
$(1,1)$-classes these conditions induce relations between 
$F$
(written in terms of a base of
$\Lambda$)
and the period matrix
$\tau$. 
In a family 
$\left\{ A_t \right\}_{t \in \Delta}$
of complex tori these relations must be satisfied for 
$\tau_t$,
$t \neq 0$,
if a holomorphic vector bundle
$E$
on
$A_0$ still has a holomorphic structure on
$A_t$.

\noindent
But the existence problem for vector bundles of higher rank with prescribed 
Chern classes is much more difficult than in the case of line bundles. On 
(non-algebraic) complex tori this problem is completely 
solved only in dimension 2 and rank 2 \cite{Tom99}, \cite{TT02}. 
Consequently, to construct
a counter-example to Kodaira's conjecture with vector bundles of higher rank 
it is not enough to give a set of Chern classes and to prove that these classes
can be Chern classes only for isolated period matrix. On the other hand if 
there is a positive family of such period matrices there may be still a 
counter-example depending on the existence of vector bundles with these
Chern classes only on isolated members of the families. 

\noindent
Finally the two simplest cases of this general setting are considered.

\subsection{Line bundles in arbitrary dimensions}

\noindent
Let
$X$
be a complex torus of dimension
$g$
given by the period matrix
$(\tau,\ \mathbf{1}_g)$.
By the characterization of the N\'eron-Severi group in the last section a 
skew symmetric matrix 
$E \in M_{2g}(\mathbb{Z})$
is a 
$(1,1)$-form iff the entries of
$\tau$
satisfy
$\left( \begin{array}{c}
        g \\
        2
        \end{array} \right)$
equations. Consequently, 3 skew symmetric matrices
$E_1, E_2, E_3 \in M_{2g}(\mathbb{Z})$
should determine at most a finite number of 
$g \times g$
period matrices
$\tau$
such that
$E_1, E_2, E_3$
are first Chern classes of line bundles on
$X_\tau$.

\noindent
As in the last section, for given
$g$
one can choose random entries for 
$E_1, E_2, E_3$
and compute the locus
$V(E_1) \cap V(E_2) \cap V(E_3)$
of
$\tau$'s
as above. But already in dimension 4 this locus turns out to be empty for
randomly chosen entries. This means that only special triples of matrices
belong to the N\'eron-Severi group of a complex torus, and it seems difficult 
to find one such that furthermore the above locus consists of isolated points.
And then one has still to prove that the period matrices in this locus 
determine a non-algebraic complex torus.

\subsection{Rank 2 vector bundles in dimension 3}

\noindent
This is the simplest case with vector bundles of rank 
$>1$. 
Unfortunately, by Poincar\'e duality
\[ H^{2,2}(X, \mathbb{Z}) \cong  H^{1,1}(X, \mathbb{Z}),
   H^{3,1}(X, \mathbb{Z}) \cong  H^{0,2}(X, \mathbb{Z}),
   H^{1,3}(X, \mathbb{Z}) \cong  H^{2,0}(X, \mathbb{Z}) \]
and the equations for a skew symmetric matrix in
$M_{2g}(\mathbb{Z})$
to be a 
$(2,2)$-form do not differ from those for
$(1,1)$-forms. Hence in this case a counter-example may be found only by 
closer 
considering the question for which complex tori exist rank 2 vector bundles with 
given Chern classes.

\vspace{0.2cm}

\noindent
Of course more difficult settings starting with rank 2 vector bundles on 
4-dimensional complex tori may give positive results. On the other hand the 
examples above give enough evidence to turn around one's point of view, into
an attempt to prove Kodaira's conjecture (at least in these special cases).

\section{Deformations of Projective Bundles}

\noindent
In this final section we generalize Proposition 2:

\begin{thm}
Let $X$ be a compact K\"ahler manifold which has a ${\mathbb{P}_r}$-bundle 
structure $X\to Y$ over some compact K\"ahler manifold $Y$. Then for every 
deformation $\mathcal{X}\to S$ with $\mathcal{X}_0\simeq X$, the nearby 
fibers $\mathcal{X}_t$ have a ${\mathbb{P}_r}$-bundle structure
$\mathcal{X}_t\to \mathcal{Y}_t$ where $\mathcal{Y}$ is a deformation of 
$Y$ in a neighborhood of $t=0$. Moreover, if $X=\mathbb{P}(V)$ for some 
vector bundle $V$ on $A$, then $\mathcal{X}_t=\mathbb{P}(V_t)$ for a 
suitable deformation $V_t\to \mathcal{A}_t$ of $V\to A$.
\end{thm}

\begin{proof} Let $q: \cC \to \cT$ be the irreducible component of the cycle space relative
to $\pi: \mathcal{X} \to S$ containing the fibers of $X \to Y.$ So $\cT$ parametrizes 
deformations of the $\mathbb{P}_r$ to nearby fibers $X_y.$ 
Since the normal bundle in $\cC$ to these projective spaces is trivial, it follows
immediately that (possibly after shrinking $S$) $\cT$ is smooth.
Let $p: \cC \to X$ denote
the projection and notice that there is another
canonical projection $r: \cT \to S$ realizing $\cT$ as a family $(T_s).$ We will also
consider $\cC_s = q^{-1}(T_s)$ with projection $q_s$ to $T_s.$ Now $q_0$ is a
$\mathbb{P}_r-$bundle. Therefore for small $s$ also the maps $q_s$ are first submersions
and second projective bundles (since projective space is locally rigid). 
Having in mind that $p_0: \cC_0 \to X_0$
is an isomorphism, we see that $p$ is an isomorphism so that
all $X_t$ are projective bundles for small $t.$
\\ The vector bundle statement finally is proved just as in Proposition 2.
\end{proof}

\end{document}